\documentclass[a4paper]{article}

\makeatletter
\RequirePackage{ifthen}
\newcommand{\Draft}{1}
	\AtBeginDocument{\@ifpackageloaded{showkeys}
	{\renewcommand{\Draft}{1}}
	{\renewcommand{\Draft}{0}}
}
\makeatother

\usepackage{amsmath}
\usepackage{amssymb}
\usepackage{amsthm}
\usepackage{bm}

\usepackage{dsfont}
\usepackage{enumerate}
\usepackage{epic}
\usepackage{eepic}
\usepackage[OT2,OT1]{fontenc}
\usepackage{ifthen}
\usepackage[latin1]{inputenc}
\usepackage{pifont}
\usepackage{rotating}
\usepackage{stmaryrd}
\usepackage{textcomp}
\usepackage{theoremref}
\usepackage{tikz}
	\usetikzlibrary{arrows}
	\usetikzlibrary{patterns}
\usepackage{todonotes}
\usepackage{url}

\newcommand{\bb}[1]{{\mathbb{#1}}}

\begingroup\makeatletter\ifx\SetFigFont\undefined%
\gdef\SetFigFont#1#2#3#4#5{%
  \reset@font\fontsize{#1}{#2pt}%
  \fontfamily{#3}\fontseries{#4}\fontshape{#5}%
  \selectfont}%
\fi\endgroup%

\newcommand\cyr{%
\renewcommand\rmdefault{wncyr}%
\renewcommand\sfdefault{wncyss}%
\renewcommand\encodingdefault{OT2}%
\normalfont
\selectfont}
\DeclareTextFontCommand{\textcyr}{\cyr}

\DeclareMathOperator{\RE}{Re}
\DeclareMathOperator{\IM}{Im}
\renewcommand{\Re}{\RE}
\renewcommand{\Im}{\IM}

\newlength{\maxlabwidth}

\theoremstyle{plain}
	\newtheorem{lemma}{Lemma}
	
	\newtheorem{theorem}{Theorem}
	
	\newtheorem{ntheoreM}[lemma]{}
\theoremstyle{definition}
	\newtheorem{definitioN}[lemma]{Definition}
\theoremstyle{remark}
	\newtheorem{remarK}[lemma]{Remark}
	\newtheorem{examplE}[lemma]{Example}
	\newtheorem{nremarK}[lemma]{}
\newcommand{\thlab}[1]{\thlabel{#1}\label{#1 }}
\renewcommand{\qedsymbol}{\raisebox{-2pt}{\large\ding{113}}}
\newcommand{\defendsymbol}{$\lozenge$}
\newcommand{\qedsymbolsave}{\qedsymbol}

\newenvironment{remark}{\begin{remarK}}{
	\renewcommand{\qedsymbolsave}{\qedsymbol}\renewcommand{\qedsymbol}{\defendsymbol}
	\popQED{\qed}\renewcommand{\qedsymbol}{\qedsymbolsave}\end{remarK}}

\begin{document}

\ifthenelse{\Draft=1}{\pagenumbering{roman}}

\begin{flushleft}
	{\Large\bf Density of the spectrum of Jacobi matrices}
	\\[1.5mm]
	{\Large\bf with power asymptotics}
	\\[4mm] 
	\textsc{Raphael Pruckner \, \
		\hspace*{-19pt}
		\renewcommand{\thefootnote}{\fnsymbol{footnote}}
		\footnote{raphael.pruckner@tuwien.ac.at}
		\setcounter{footnote}{0}
	}
	\\
	{\footnotesize
	\vspace*{4mm}
	Institute for Analysis and Scientific Computing, Vienna University of Technology\\
	Wiedner Hauptstra{\ss}e\ 8--10/101, 1040 Wien, AUSTRIA
	}
	\end{flushleft}
	\vspace*{0mm}
	{\small \noindent
	\textbf{Abstract:}
		We consider Jacobi matrices $J$ 
		whose parameters have the power asymptotics 
		$\rho_n=n^{\beta_1} \left(  x_0 + \frac{x_1}{n} + {\rm O}(n^{-1-\epsilon})\right)$
		and 
		$q_n=n^{\beta_2} \left(  y_0 + \frac{y_1}{n} + {\rm O}(n^{-1-\epsilon})\right)$ for the off-diagonal and 
		diagonal, respectively.
		
		We show that for $\beta_1 > \beta_2$, 
		or $\beta_1=\beta_2$ and $2x_0 > |y_0|$, 
		the matrix $J$ is in the limit circle case and 
		the convergence exponent of its spectrum is $1/\beta_1$.
		Moreover, we obtain upper and lower bounds for the upper density of the spectrum.
		
		When the parameters of the matrix $J$ have
		a power asymptotic with one more term,
		we characterise the occurrence of the limit circle case completely 
		(including the exceptional case $\lim_{n\to \infty} |q_n|\big/ \rho_n = 2$)
		and determine the convergence exponent in almost all cases.
 	\\[3mm] 
	}
	\\[1mm]

\ifthenelse{\Draft=1}{
\fbox{
\parbox{100mm}{
\hspace*{0pt}\\
\centerline{{\Large\ding{45}}\quad\,{\large\sc Draft}\quad{\Large\ding{45}}}
\hspace*{0pt}\\[-3mm]
In final version:\qquad \ding{233} disable package \textsf{showkeys}.\\
Download latest version of bibtex-database-woracek.bib:\\
\centerline{\textsf{asc.tuwien.ac.at/\~{}funkana/woracek/bibtex-database-woracek.bib}}
\\[2mm]
\textcircledP\ Preliminary version Thu 27 Sep 2018 07:22
\\[2mm]
\ding{229}\quad to obtain proper eepic-pictures use latex/dvipdf to compile
\\[-2mm]
}
}
\tableofcontents
\listoftodos
%
\newpage
\pagenumbering{arabic}
\setcounter{page}{1}
}
{}


%
%
%
\section{Introduction}
A Jacobi matrix $J$ is a tridiagonal semi-infinite matrix
\begin{equation*} 
	J=
	\begin{pmatrix}
		q_0 & \rho_0 && \\
		\rho_0 & q_1 & \rho_1 &  \\ 
		& \rho_1 & q_2 & \smash{\ddots} & \\
		&& \ddots & \ddots 
	\end{pmatrix}
\end{equation*}
with real $q_n$ and positive $\rho_n$. A Jacobi matrix induces a closed symmetric operator $T_J$ on $\ell^2(\bb N)$, 
namely as the closure of the natural action of $J$ on the subspace of finitely supported sequences, 
see, e.g.,\ \cite[Chapter~4.1]{akhiezer:1961}. 
There occurs an alternative: 
\begin{itemize}
	\item $T_J$ is selfadjoint; one speaks of the limit point case (lpc), or, in the language of \cite{akhiezer:1961}, type D.
	\item $T_J$ has defect index $(1,1)$ and is entire in the sense of M.G.Kre{\u\i}n; one speaks of the limit circle case (lcc), or, 
	synonymously, type C.
\end{itemize}
In the lpc, the spectrum of $T_J$ may be discrete, continuous, or be composed of different types.
If $J$ is in the lcc, then the spectrum of every canonical selfadjoint extension of $T_J$ is discrete, 
and each two spectra are interlacing.
In this case, we fix one such extension and denote its spectrum by $\sigma(J)$.

In general it is difficult to decide from the parameters $\rho_n,q_n$ whether $J$ is in the lcc or lpc. 
Two classical necessary conditions for occurrence of lcc are 
Carleman's condition which says that $\sum_{n=0}^\infty\rho_n^{-1}=\infty$ implies lpc, 
cf.\ \cite{carleman:1926}, and Wouk's theorem that a dominating diagonal in the sense that 
$\sup_{n\geq 0}(\rho_n+\rho_{n-1}-q_n)<\infty$ 
or $\sup_{n\geq 0}(\rho_n+\rho_{n-1}+q_n)<\infty$ implies lpc, cf.\ \cite{wouk:1953}. 
A more subtle result, which gives a sufficient condition for lcc, is due to Yu.M.Berezanski{\u\i}, 
cf.\ \cite[Theorem~4.3]{berezanskii:1956}, \cite[VII,Theorem~1.5]{berezanskii:1968} or \cite[Addenda~5.,p.26]{akhiezer:1961}: 
Assume 
that $\sum_{n=0}^\infty\rho_n^{-1}<\infty$, 
that the sequence  $(q_n)_{n=0}^\infty$ of diagonal parameters is bounded,
and that the sequence $(\rho_n)_{n=0}^\infty$ of off-diagonal parameters behaves regularly in the sense that $\rho_n^2\geq\rho_{n+1}\rho_{n-1}$ (log-concavity). 
Then $J$ is in the lcc.
An extension and modern formulation of this result can be found in \cite[Theorem~4.2]{berg.szwarc:2014}.
In particular, instead of $(q_n)_{n=0}^\infty$ being bounded it is enough to require $\sum_{n=0}^\infty |q_n|\big/\rho_n<\infty$.

There is a vast literature dealing with Jacobi matrices in the lpc, whose aim is to establish discreteness of the spectrum 
and investigate spectral asymptotics, e.g.,\ \cite{boutet.zielinski:2012,deift:1999,janas.malejki:2007,janas.naboko:2004,tur:2003}.
Contrasting this, if $J$ is in the lcc, not much is known about the asymptotic behaviour of the spectrum. 

The probably first result in this direction is due to M.Riesz\ \cite{riesz:1923a} and states that
the spectrum of a Jacobi matrix in the lcc case
is sparse compared to the integers, in the sense that
($\lambda_n^\pm$ denote the sequences of positive or negative points in $\sigma(J)$, arranged according to increasing modulus)
\[
	\lim_{n\to\infty}\frac n{\lambda_n^\pm}=0
	.
\]
A deeper result 
holds in the context of
the already mentioned work of Berezanski{\u\i}.
Under the mentioned assumptions,
C.Berg and R.Szwarc showed 
that the convergence exponent $\rho(\sigma(J))$ of the spectrum coincides with 
the convergence exponent of the sequence $(\rho_n)_{n=0}^\infty$ of off-diagonal parameters of $J$,
cf.\ \cite[Theorem~4.11]{berg.szwarc:2014}.
Recall that the convergence exponent of a sequence $(x_n)_{n=0}^\infty$, which we denote by $\rho((x_n)_{n=0}^\infty)$,
is defined as  the greatest lower bound of all $\alpha>0$ such that $\sum_{n=0}^\infty x_n^{-\alpha} < \infty$.

In this paper we contribute to the study of the spectrum of Jacobi matrices in the lcc. 
We investigate Jacobi matrices $J$ whose parameters have, for some $\epsilon>0$, power asymptotics
\begin{equation} \label{P32}
	\rho_n = n^{\beta_1} \Big(  x_0 + \frac{x_1}{n} + {\rm O}\Big(\frac{1}{n^{1+\epsilon}}\Big)\Big)
	, \quad
	q_n = n^{\beta_2} \Big(  y_0 + \frac{y_1}{n} +  {\rm O}\Big(\frac{1}{n^{1+\epsilon}}\Big)\Big)
	,
\end{equation} 
with $x_0>0$, $y_0\neq 0$. 

%

In our first theorem, 
we show that,
apart from the exceptional case that $\lim_{n\to \infty} |q_n|\big/\rho_n =2$,
a characterisation of the lcc is possible.
Moreover, we give bounds for the upper density of the spectrum, in particular, determine the convergence exponent.

Our second theorem treats the exceptional case that $\lim_{n\to \infty} |q_n|/\rho_n =2$,
or equivalently $\beta_1=\beta_2$ and $2x_0 = |y_0|$. 
Under a stronger assumption, we fully characterise the occurrence of the lcc and determine 
the convergence exponent of the spectrum in almost all cases .

\begin{theorem} \thlab{P1} 
Let $J$ be a Jacobi matrix 
with off-diagonal $(\rho_n)_{n=0}^\infty$ and diagonal $(q_n)_{n=0}^\infty$
which have the power asymptotics \eqref{P32}
with some $\epsilon>0$, $\beta_1,\beta_2 \in \bb R$, $x_0>0$, $y_0\neq 0$ and $x_1,y_1 \in \bb R$.

Consider the following two cases.
\begin{enumerate}[$(i)$]
\item $\beta_1 \leq \beta_2$, and $2x_0<|y_0|$ if $\beta_1=\beta_2$

In this case, $J$ is in the lpc. 
\item $\beta_1 \geq \beta_2$, and $2x_0>|y_0|$ if $\beta_1=\beta_2$

In this case, $J$ is in the lcc if and only if $\beta_1$ is greater than $1$. 
\end{enumerate}
In the lcc, the convergence exponent of the spectrum is $1/\beta_1$.
Moreover, we have the following bounds for the upper density of the spectrum,
\begin{equation*} 
	\frac{\beta_1 - 1}{\beta_1}  
	\left(\frac{1}{x_0}\right)^{1/\beta_1}	
	\leq 
	\limsup_{r\to \infty} \frac{n_{\sigma} (r)}{r^{1/{\beta_1}}} 
	\leq 
	e\frac{\beta_1 }{\beta_1 - 1} \left( \frac{a}{x_0}\right)^{1/\beta_1}
	,
\end{equation*}
where $n_\sigma(r):=\# (\sigma(J) \cap [-r,r])$ denotes the counting function, and
\[
	a:=\begin{cases}
	1 &, \beta_1>\beta_2, \\
	 \Big(1-\frac{y_0^2}{4 x_0^2}\Big)^{-1/2} \ \  &, \beta_1=\beta_2
	.
	\end{cases}
\]
\end{theorem}

\begin{remark}
Note that case $(i)$ in \thref{P1} is equivalent to $\lim_{n\to \infty} |q_n|/\rho_n >2$,
whereas case $(ii)$ corresponds to $\lim_{n\to \infty} |q_n|/\rho_n <2$. 
\end{remark}

\begin{remark}
Having \eqref{P32} implies that $(\rho_n)_{n=0}^\infty$ is log-concave. Hence, if $\beta_1>\beta_2+1$, 
the above discussed extension of Berezanski{\u\i}'s theorem applies and yields $\rho(\sigma(J))=\rho((\rho_n)_{n=0}^\infty)$.
Note that the convergence exponent of a sequence $(\rho_n)_{n=0}^\infty$ with \eqref{P32} is $1/\beta_1$.

For $\beta_1 > \beta_2+1$, Theorem 1 
refines this result by providing explicit estimates for 
the upper density of the spectrum.
However, the main significance is that the statement remains valid for $\beta_1>\beta_2$, 
and even in some cases where $\beta_1=\beta_2$, i.e., where diagonal and off-diagonal parameters are comparable.
\end{remark}

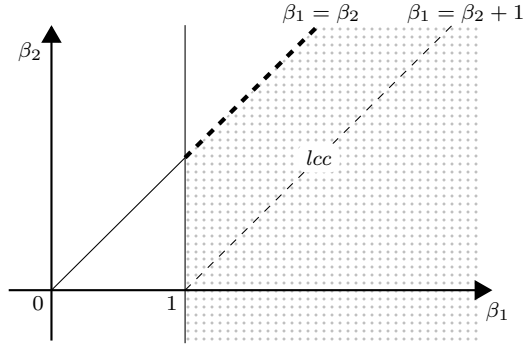
\begin{center}
	{ \footnotesize
		%
		\begin{tikzpicture}[x={(20pt,0pt)},y={(0pt,20pt)},scale=0.5]
		\draw[-triangle 60,thick] (2,0.1) -- (2,12);				
		\draw[-triangle 60,thick] (0.4,2) -- (18.5,2);
		\draw (1.5,1.5) node {$0$};
		\draw (6.5,1.5) node {$1$};
		\draw (18.7,1.2) node {$\beta_1$};
		\draw (1.2,11) node {$\beta_2$};
		\draw (7,0) -- (7,12);
		
		\draw (2,2) -- (7,7);
		\draw[dashed,ultra thick] (7,7) -- (12,12);	
		\draw[dashed] (7,2) -- (17,12);			
		\fill[pattern=dots,opacity=0.5]
		(7,0) -- (7,7) -- (12,12) -- (18,12) -- (18,0);
		
		\draw (12.1,12.45) node {$\beta_1=\beta_2$};
		\draw (17.5,12.45) node {$\beta_1=\beta_2+1$};
		\draw (12,7) node[fill=white, rounded corners=4pt] {\emph{lcc}};	
		\end{tikzpicture}
		%
	}
\end{center}

In order to handle the exceptional case, we require
the stronger assumption that the parameters of the Jacobi matrix $J$ have, for some $\epsilon>0$, power asymptotics of the form
\begin{equation} \label{P32E}
	\rho_n = n^{\beta_1} \Big( x_0 + \frac{x_1}{n} + \frac{x_2}{n^2}+ {\rm O}\Big(\frac{1}{n^{2+\epsilon}}\Big)\Big)
	,\quad
	q_n = n^{\beta_2} \Big(  y_0 + \frac{y_1}{n} + \frac{y_2}{n^2} + {\rm O}\Big(\frac{1}{n^{2+\epsilon}}\Big)\Big)
	,
\end{equation}
with $x_0>0$ and $y_0\neq 0$.

\begin{theorem} \thlab{ADD1}
Let $J$ be a Jacobi matrix 
with off-diagonal $(\rho_n)_{n=0}^\infty$ and diagonal $(q_n)_{n=0}^\infty$
which have the power asymptotics \eqref{P32E}
with some $\epsilon>0$, $\beta_1,\beta_2 \in \bb R$, $x_0>0$, $y_0\neq 0$ and $x_1,y_1,x_2,y_2 \in \bb R$.

Assume that $\beta_1=\beta_2=:\beta$ and $2x_0 = |y_0|$.
Then, exactly one of the following cases takes place.
\begin{enumerate}[$(i)$]
	\item ${\beta\in \big(-\infty,\frac 3 2\big] \cup \big(\frac{2x_1}{x_0}-\frac{2y_1}{y_0}, \infty\big)}$
	
	In this case, $J$ is in the lpc.

	\item $\beta\in \big(\tfrac 3 2, \frac{2x_1}{x_0}-\frac{2y_1}{y_0}\big)$
	
	In this case, $J$ is in the lcc. Regarding the convergence exponent of the spectrum, we have
	\begin{equation} \label{P27}
	\rho(\sigma(J))
		\begin{cases}
			\in [\frac{1}{\beta},\frac{1}{2(\beta-1)}]   &, \frac{3}{2}< \beta < 2, \\
			= \frac{1}{\beta} &, \beta \geq 2.
		\end{cases}
	\end{equation}
	\item $\beta = \frac{2x_1}{x_0}-\frac{2y_1}{y_0}>\tfrac 3 2 $
	
	In this case, $J$ is in the lcc if and only if 
	$2 < \beta < \frac{3}{2} + \frac{2z_2}{x_0}$,
	where
	\[
		z_2=x_0 \left(
		\frac{2x_2}{x_0}-\frac{2y_2}{y_0}+\frac{y_1}{y_0} -\frac{2x_1 y_1}{x_0 y_0} + \frac{2 y_1^2}{y_0^2}
		\right)
		.
	\]
	In the lcc, the convergence exponent of the spectrum is $1/\beta$.
\end{enumerate}
When $J$ is in the lcc, 
the following lower estimate of the density
of the spectrum holds,
\begin{equation*}
	\frac{\beta - 1}{\beta}  
	\left(\frac{1}{x_0}\right)^{1/\beta}	
	\leq \limsup_{r\to \infty} \frac{n_{\sigma} (r)}{r^{1/{\beta}}}
	.
\end{equation*}
\end{theorem}

\begin{remark}
We strongly believe that the convergence exponent of the spectrum is equal to $1/\beta$ whenever $J$ is in the lcc, even for $\beta\in (\tfrac 3 2, 2)$ in the case $(ii)$.
\end{remark}

\begin{remark}
In case $(iii)$ of \thref{ADD1}, the parameter $\beta$ is already given by  $x_0,y_0,x_1,y_1$.
Hence, the condition
$2 < \beta < \frac{3}{2} + \frac{2z_2}{x_0}$
can equivalently be written as the two conditions
\[
	 1 < \frac{x_1}{x_0}-\frac{y_1}{y_0}, \quad \quad \frac{x_1}{|y_0|} + \frac{y_1}{y_0}\left( \left(\frac{x_1}{x_0}-\frac{y_1}{y_0}\right) -1 \right) <\frac{3}{8} + 	\frac{x_2}{x_0} - \frac{y_2}{y_0}
	.
\]
Moreover, the notation $z_2$, which is introduced in that case, relates to Wouk's theorem, cf.\ \eqref{P2}.
\end{remark}

%

In the proofs of these theorems we first establish that the power asymptotics of the Jacobi parameters, i.e.\ \eqref{P32} or \eqref{P32E},  give rise to 
the asymptotic behaviour of 
a fundamental solution of the finite difference equation \eqref{P3} corresponding to $J$.
This is achieved by applying theorems of R.Kooman.
In the proof of \thref{P1}, we use \cite[Corollary~1.6]{kooman:1998}, 
which is a generalisation of the classical Poincare-Perron theorem to the case that the zeros of the characteristic equation may have the same modulus but are distinct. 
In the exceptional case, the characteristic equation has a double zero and the more involved theorem \cite[Theorem~1]{kooman:2007} is needed.
In any case, the asymptotic behaviour of solutions directly leads to a characterisation of the lcc.

The crucial step is to estimate the upper density of the spectrum and determine the convergence exponent in the lcc.
Here we use the fact that the growth of the counting function $n_\sigma$ 
relates to the growth of the canonical product having the spectrum as its zero-set. 
The upper density of the spectrum is in our setting always bounded from below 
by a result of C.Berg and R.Szwarc, i.e.\ \cite[Proposition~7.1]{berg.szwarc:2014}.
In the proof of \thref{P1},
we obtain an upper bound of the upper density by estimating the canonical product by hand.
In particular, this determines the convergence exponent of the spectrum.
If one is only interested in the convergence exponent, it is enough to 
apply \cite[Theorem~1.2]{berg.szwarc:2014}.
In the situation of \thref{ADD1}, both approaches fail and a better estimate of the canonical product is needed, cf.\ \thref{P28}.
This is achieved by writing the Jacobi matrix as a Hamburger Hamiltonian of a canonical system
and applying \cite[Theorem~2.7]{pruckner.romanov.woracek:2016}, 
which goes back to a theorem of R.Romanov, cf.\ \cite[Theorem~1]{romanov:2017}.
\section{Proof of Theorem~1}
%
%

%



Let $J$ be a Jacobi matrix 
with parameters $\rho_n$ and $q_n$ 
having the power asymptotics \eqref{P32}
with some $\epsilon>0$, $\beta_1,\beta_2 \in \bb R$, $x_0>0$, $y_0\neq 0$ and $x_1,y_1 \in \bb R$.

Recall Wouk's theorem, which is formulated in the Introduction. Since $q_n$ does not change its sign for $n$ large enough, this theorem states that if
$\sup_{n\geq 1} (\rho_n + \rho_{n-1} - |q_n|)<\infty$, then $J$ is in the lpc. 

In case $(i)$, we have $\lim_{n\to \infty} |q_n|/\rho_n >2$, which implies
\[
	\lim_{n\to \infty} \frac{\rho_n + \rho_{n-1} - |q_n|}{\rho_n} = 2-\lim_{n\to \infty} \frac{|q_n|}{\rho_n}  < 0
	.
\]
Hence, $J$ is in the lpc by Wouk's theorem.
It remains to treat case $(ii)$. 
Thus, assume $\delta:=\beta_1-\beta_2\geq 0$, and $2x_0 > |y_0|$ if $\delta=0$.

\subsubsection*{\underline{Step 1:} Growth of solutions}
We start with studying asymptotics of solutions of the difference equation
\begin{equation} \label{P3}
\rho_{n+1} u_{n+2} + q_{n+1} u_{n+1} + \rho_n u_n = 0
.
\end{equation}
Dividing by $\rho_{n+1}$ yields
\[
u_{n+2} + C_1(n) u_{n+1} + C_0(n) u_n = 0
,
\]
with
\begin{align*}
	&C_1(n)=\frac{q_{n+1}}{\rho_{n+1}} = n^{-\delta} 
	\Big( \frac{y_0}{x_0}   + \frac{1}{n}\Big(\frac{y_1-\delta y_0}{x_0} - \frac{y_0 x_1}{x_0^2}\Big) + {\rm O}\Big(\frac{1}{n^{1+\epsilon}}\Big) \Big)
	,\\
	&C_0(n)=\frac{\rho_n}{\rho_{n+1}} = 1 - \frac{\beta_1}{n} + {\rm O}\Big(\frac{1}{n^{1+\epsilon}}\Big)
	.
\end{align*}
We denote by $\alpha_1(n), \alpha_2(n)$ the zeros of the characteristic polynomials, i.e.
\[
	 x^2+C_1(n)x+C_0(n)=0
	.
\]
Note that the limit
\[
\lim_{n\to \infty} \frac{C_1(n)^2}{4} - C_0(n) 
=
\begin{cases}
-1 \ &, \delta >0,\\
\frac{y_0^2}{4x_0^2}-1 &, \delta=0,
\end{cases}
\]
is negative by assumption. 
It follows that $\alpha_1(n)$ and $\alpha_2(n)$ are, for $n$ large enough, complex conjugate numbers,
which converge to distinct numbers, i.e.
\begin{equation} \label{P12}
	\lim_{n\to \infty} \alpha_{1,2}(n) = \begin{cases}
	\pm i &, \delta >0,\\
	\frac{-y_0}{2 x_0} \pm i \sqrt{1-\frac{y_0^2}{4 x_0^2}} &, \delta=0.
	\end{cases}
\end{equation}
Moreover, $(C_i(n)-C_i(n+1))$ is summable for $i=0,1$, due to
\begin{align*}
	|C_1(n)-C_1(n+1)| &\leq 
		C \big| n^{-\delta} - (n+1)^{-\delta}\big| + {\rm O}(n^{-1-\epsilon}) 
		\\
		&\leq	C \delta n^{-\delta -1} + {\rm O}(n^{-1-\epsilon}), \\
	|C_0(n)-C_0(n+1)| &= \Big|-\frac{\beta_1}{n}  + \frac{\beta_1}{n+1} \Big| + {\rm O}(n^{-1-\epsilon}) = {\rm O}(n^{-1-\epsilon})
	.
\end{align*}
Now \cite[Corollary~1.6]{kooman:1998}
yields two linearly independent, complex conjugate solutions $(v_n^{(1)})_{n=1}^\infty$, $(v_n^{(2)})_{n=1}^\infty$  
of \eqref{P3} with
\begin{equation} \label{P19}
	v_n^{(i)}= (1+{\rm o}(1)) \prod_{k=0}^{n-1} \alpha_i(k),
	\quad n\in \bb N
	.
\end{equation}
By using \cite[Lemma~4]{kooman:2007}, adding a summable perturbation, we get
\begin{equation} \label{P15}
	\Big|\prod_{k=0}^{n-1} \alpha_i(k) \Big|^2
	=\prod_{k=0}^{n-1} C_0(k) 
	=\prod_{k=0}^{n-1} \Big( 1 - \frac{\beta_1}{k} + {\rm O}\Big(\frac{1}{k^{1+\epsilon}}\Big)\Big)
	= (d^2+{\rm o}(1)) n^{-\beta_1}
	,
\end{equation}
for a constant $d>0$. 
Hence, the normalized solutions $u_n^{(j)}:=v_n^{(j)}/d$ for $j=1,2$ satisfy $|u_n^{(1)}|=|u_n^{(2)}|= (1+{\rm o}(1)) n^{-\beta_1/2}$. 
In particular, 
they are square-summable and $J$ is in the lcc if and only if $\beta_1>1$.

\subsubsection*{\underline{Step 2:} The lower bound in the lcc}

From the first step we know that the corresponding moment problem is in the lcc. Thus, the  Nevanlinna matrix $\big(\begin{smallmatrix} A(z) & B(z) \\ C(z) & D(z) \end{smallmatrix} \big)$ which parametrizes all solutions of the moment problem is available, cf.~\cite{nevanlinna:1922, akhiezer:1961}.
This four entries are canonical products and have the same growth, i.e.\ the same type w.r.t.\ any growth function, cf.~\cite[Proposition~2.3]{baranov.woracek:smsub}.
In particular, they have the same order and type.

The zeros of $B$ interlace with the spectrum of any canonical selfadjoint extension of $T_J$. Thus, the counting function of the zeros of $B$, which we denote by $n_{B}(r)$, differs from $n_{\sigma} (r)$ by at most $1$. 
Hence, knowledge about the growth of any entry of the Nevanlinna matrix can be used to derive knowledge about the distribution of the spectrum.
In particular, the order of $B$ coincides with the convergence exponent $\rho(\sigma(J))$ of the spectrum, and the type of $B$ is comparable 
to the upper density of the spectrum with explicit constants.

The order and type of the entries of the Nevanlinna matrix are, by \cite[Proposition~7.1~(\textit{ii}),(\textit{iii})]{berg.szwarc:2014}, bounded from below by the order and type of the entire function
\[
	H(z)= \sum_{n=0}^\infty b_{n,n}z^n
	,
\]
where $b_{n,n}=\left( \rho_1 \rho_2 \cdot \ldots \cdot \rho_{n-1} \right)^{-1}$ denotes 
the leading coefficient of the $n$-th orthogonal polynomial of the first kind, denoted by $P_n(z)$.
The power asymptotic of $\rho_n$ yields
\[
	b_{n,n} = 
	(C+{\rm o}(1))\big[n!\big]^{-\beta_1} n^{\beta_1 - \frac{x_1}{x_0}} x_0^{-n+1}
	,
\]
for a constant $C>0$. 
By the standard formula for the order and type of a power series, cf.\ \cite[Theorem~2]{levin:1980}, we get
that the order of $H(z)$ is $1/{\beta_1}$, and the type w.r.t.\ this order 
is equal to $\beta_1 x_0^{-1/\beta_1}$. 

Thus, we get $\rho(\sigma(J))\geq 1/\beta_1$ and $\tau(B) \geq \beta_1 x_0^{-1/\beta_1}$,
where $\tau(B)$ denotes the type of $B$ w.r.t.\ the order $1/\beta_1$. 
The inequality between the type of a canonical product
and the upper density of its zeros, 
cf.~\cite[eq.\ (1.25)]{levin:1980}, gives
\[
\limsup_{r\to \infty} \frac{n_{\sigma} (r)}{r^{1/{\beta_1}}}
=
\limsup_{r\to \infty} \frac{n_{B} (r)}{r^{1/{\beta_1}}} 
 \geq \frac{\beta_1 - 1}{\beta_1^2}	\, \tau(B)  \geq  \frac{\beta_1 - 1}{\beta_1} 
\Big(\frac{1}{x_0}\Big)^{1/\beta_1}
 .
\]

\subsubsection*{\underline{Step 3:} The upper bound in the lcc}

In the first step we have seen that the difference equation \eqref{P3} has a fundamental system of solutions
with $|u_n^{(j)}| = (1+{\rm o}(1)) n^{-\beta_1/2}$ for $j=1,2$.

The orthogonal polynomials of the first and second kind associated with the matrix $J$, denoted by $P_n(0)$ and $Q_n(0)$ respectively, 
are also linearly independent solutions of \eqref{P3}. 
Therefore, both $\big(P_n^2(0)\big)_{n=1}^\infty$ and $\big(Q_n^2(0)\big)_{n=1}^\infty$ 
are in $l^p$ for $p>1/\beta_1$.
By \cite[Theorem~1.2]{berg.szwarc:2014} the order of $B$ is at most, and hence equal to, $1/\beta_1$.

We are going to estimate the density of the spectrum from above by analysing the growth of $B$ more precisely. 
To this end, we write the Nevanlinna matrix as a product, i.e.
\[
	\begin{pmatrix}	A_{n+1}(z) & B_{n+1}(z) \\ C_{n+1}(z) & D_{n+1}(z) 	\end{pmatrix}
	= 
	\big( I + z 
	R_n
	\big)
	\begin{pmatrix}	A_{n}(z) & B_{n}(z) \\ C_{n}(z) & D_{n}(z) 	\end{pmatrix}
	,
\]
with
\[
	R_n:=\begin{pmatrix}-P_n(0)Q_n(0) & Q_n^2(0) \\ -P_n^2(0) & P_n(0)Q_n(0)	\end{pmatrix}
	=	\begin{pmatrix}
		Q_n(0) \\ P_n(0)
	\end{pmatrix} 	\begin{pmatrix}
		Q_n(0) \\ P_n(0)
	\end{pmatrix}^T
	\begin{pmatrix}0 & 1 \\ -1 & 0 \end{pmatrix}
	.
\]
Here, $A_n, B_n, C_n$ and $D_n$ are polynomials,
which converge to the corresponding entry of the Nevanlinna matrix, cf.\ \cite[p.~14/54]{akhiezer:1961} and \cite[eq.\ (39)]{berg.szwarc:2014}.
Hence, 
the spectral norm of the Nevanlinna matrix can be written as
\begin{equation}\label{P33}
	\left\|\begin{pmatrix}	A(z) & B(z) \\ C(z) & D(z) 	\end{pmatrix}\right\|
	= \bigg\| \prod_{n=0}^{\infty} (I + z R_n) \bigg\| 
	.
\end{equation}
Let $T$ denote the regular $2\times 2$ matrix such that
\[
	\begin{pmatrix}
		Q_n(0) \\ P_n(0)
	\end{pmatrix} = T \bigg(\begin{matrix}
		u_n^{(1)} \\ u_n^{(2)}
	\end{matrix} \bigg)
	.
\]
Before we use the submultiplicativity of the norm on the right-hand side of \eqref{P33}, we rewrite the factors as follows.
\begin{align*}
I+z R_n &= I + z T 
\bigg(\begin{matrix}	u_n^{(1)} \\ u_n^{(2)}	\end{matrix} \bigg)
\bigg(\begin{matrix}	u_n^{(1)} \\ u_n^{(2)}	\end{matrix}\bigg)^{\! T}
 T^T \begin{pmatrix}	0 & 1 \\ - 1& 0 \end{pmatrix}
\\
&=T \bigg[ T^{-1} \begin{pmatrix}0 & -1 \\ 1 & 0\end{pmatrix} (T^T)^{-1} + z  
\bigg(\begin{matrix}	u_n^{(1)} \\ u_n^{(2)}	\end{matrix} \bigg)
\bigg(\begin{matrix}	u_n^{(1)} \\ u_n^{(2)}	\end{matrix}\bigg)^{\! T}
 \bigg]T^T \begin{pmatrix}	0 & 1 \\ - 1& 0 \end{pmatrix}\\
&=T \bigg[ \frac{1}{\det T} \begin{pmatrix}0 & -1 \\ 1 & 0\end{pmatrix} + z  
\bigg(\begin{matrix}	u_n^{(1)} \\ u_n^{(2)}	\end{matrix} \bigg)
\bigg(\begin{matrix}	u_n^{(1)} \\ u_n^{(2)}	\end{matrix}\bigg)^{\! T}
 \bigg]T^T \begin{pmatrix}	0 & 1 \\ - 1& 0 \end{pmatrix}\\
&= \frac{T}{\det T}\bigg[ \begin{pmatrix}0 & -1 \\ 1 & 0\end{pmatrix} + z   \det T
\bigg(\begin{matrix}	u_n^{(1)} \\ u_n^{(2)}	\end{matrix} \bigg)
\bigg(\begin{matrix}	u_n^{(1)} \\ u_n^{(2)}	\end{matrix}\bigg)^{\!T}
 \bigg]T^T \begin{pmatrix}	0 & 1 \\ - 1& 0 \end{pmatrix}
\end{align*}
When taking the product, the terms outside of the brackets give 
\[
(\det T)^{-1} T^T \begin{pmatrix}	0 & 1 \\ - 1& 0 \end{pmatrix} T =  \begin{pmatrix}	0 & 1 \\ - 1& 0 \end{pmatrix}
,
\]
which is a unitary matrix, whose spectral norm is $1$.
Also note that 
\[
\left\| \bigg(\begin{matrix}	u_n^{(1)} \\ u_n^{(2)}	\end{matrix} \bigg)
	\bigg(\begin{matrix}	u_n^{(1)} \\ u_n^{(2)}	\end{matrix}\bigg)^{\! T} \right\| = \big|	u_n^{(1)} \big|^2 +  \big|	u_n^{(2)}\big|^2
	=2 \big|	u_n^{(1)} \big|^2 
	.
\]
Hence, pulling the norm into the product in \eqref{P33} yields,
with the notation
\[
	F(z):=\prod_{n=0}^\infty (1+ z f_n n^{-\beta_1}), \quad f_n:=2 |\det T| | u_n^{(1)}|^2 n^{\beta_1}
	,
\]
nothing but
\begin{equation}  \label{P34}
	\left\|\begin{pmatrix}	A(z) & B(z) \\ C(z) & D(z) 	\end{pmatrix}\right\| 
\leq c \, F(|z|)
,
\end{equation}
for a constant $c>0$ which depends on $T$ only.
Therefore, order and type of the entries of the Nevanlinna matrix do not exceed the order and type of $F$.

Due to the first step, we have $f_n=2   |\det T| + {\rm o}(1)$.
Next, we compute the determinate of $T$ by considering the relation $g_n$
\[
	\begin{pmatrix}
		Q_{n+1}(0) & Q_{n}(0) \\ P_{n+1}(0) & P_{n}(0)
	\end{pmatrix} = T \bigg(\begin{matrix}
		u_{n+1}^{(1)}  & u_{n}^{(1)} \\ u_{n+1}^{(2)} & u_{n}^{(2)}
	\end{matrix} \bigg)
	.
\]
Taking the determinants on both sides and multiplying by $n^{\beta_1}$ gives,
due to $Q_{n+1}(0) P_{n}(0) -  P_{n+1}(0) Q_{n}(0)=\rho_{n}^{-1}$, 
\[
	\frac{n^{\beta_1}}{\rho_{n}}
	=  n^{\beta_1} \big(u_{n+1}^{(1)} u_{n}^{(2)} - u_{n+1}^{(2)}u_{n}^{(1)}\big)\det T
	.
\]
The left-hand side converges to $1/x_0$ by assumption. 
By introducing the notation $h_n:= n^{\beta_1} \big(u_{n+1}^{(1)} u_{n}^{(2)} - u_{n+1}^{(2)}u_{n}^{(1)}\big)$, we have $\det T = 1/(x_0 \lim_{n\to \infty} h_n)$.
Recall that $u_n^{(2)}$ is the complex conjugate of $u_n^{(1)}$, which gives
\begin{align} \nonumber
h_n &=
	n^{\beta_1} \Big(u_{n+1}^{(1)} \overline{u_{n}^{(1)}} - \overline{u_{n+1}^{(1)}} u_{n}^{(1)}\Big)
	= 2 i n^{\beta_1} \Im\Big(u_{n+1}^{(1)} \overline{u_{n}^{(1)}} \Big) 
	\label{P35} 
\end{align}
By \eqref{P19} and \eqref{P15} we have
\begin{align*}
	\Im\Big(u_{n+1}^{(1)} \overline{u_{n}^{(1)}}\Big)&=
	\Im\Big( (1/d^2+{\rm o }(1)) \alpha_1(n) \prod_{k=0}^{n-1} | \alpha_1(k)|^2\Big)
	\\&=(1+{\rm o }(1))n^{-\beta_1} \big(\Im(\alpha_1(n)) + {\rm o }(1) \big)
	,
\end{align*}
which gives,  due to \eqref{P12},
\[
	\lim_{n\to \infty} h_n=  2 i \lim_{n\to \infty}  \Im(\alpha_1(n)) = 
	2i /a
	,
\]
where $a$ is defined in the formulation of this theorem.
Hence, we get $\det T = a/(2 x_0 i)$ and, thus, 
\[
	f_n=	2   |\det T| + {\rm o}(1)
	=	\frac{a}{x_0} + o(1)
	.
\]
The zeros of $F$ ordered by increasing modulus behave like $(-x_0/a + o(1)) n^{\beta_1}$.
Thus, the convergence exponent of the zeros as well as the order of $F$ is equal to $1/\beta_1$.
A straight forward calculation shows that the upper density of the zeros of $F$ is equal to
$(a/x_0)^{1/\beta_1}$.
By \eqref{P34} and \cite[eq.\ (1.25)]{levin:1980}, we get the following upper bound for the type of $B$,
\[
\tau_B\leq \tau_F \leq \frac{\beta_1^2 }{\beta_1 - 1} \left( \frac{a}{x_0}\right)^{1/\beta_1}
.
\]
The fact that the type of a canonical product is, up to a constant, not lower than the 
upper density of its zeros, cf.\ \cite[Theorem~2.5.13]{boas:1954}, yields
\[
\limsup_{r\to \infty} \frac{n_{\sigma} (r)}{r^{1/{\beta_1}}}
=
\limsup_{r\to \infty} \frac{n_{B} (r)}{r^{1/{\beta_1}}} 
\leq \frac{e}{\beta_1} \tau_B \leq 
 \frac{e \beta_1 }{\beta_1 - 1} \left( \frac{a}{x_0}\right)^{1/\beta_1}
 .
\]

\vspace*{-4mm} \qed
\vspace*{-3mm}

\section{Proof of Theorem 2}

Let $J$ be a Jacobi matrix 
with parameters 
$\rho_n$ and $q_n$ having the power asymptotics \eqref{P32E}
with some $\epsilon>0$, $\beta_1,\beta_2 \in \bb R$, $x_0>0$, $y_0\neq 0$ and $x_1,y_1,x_2,y_2 \in \bb R$.
Assume that $\beta_1=\beta_2=:\beta$ and $2x_0 = |y_0|$.

A calculation shows that the expression in 
Wouk's theorem has the following power asymptotic (see also the beginning of the proof of \thref{P1}),
\begin{equation} \label{P2}
	\rho_n + \rho_{n-1} - |q_n| = n^\beta \Big(\frac{z_1}{n} + \frac{z_2}{n^2}+ {\rm O}\Big(\frac{1}{n^{2+\epsilon}}\Big) \Big),
\end{equation}
with
\begin{equation} \label{P26}
z_1=x_0\left(\frac{2x_1}{x_0} - \frac{2y_1}{y_0}   - \beta \right)
, \ \ \:
z_2=x_0 \left(
 \frac{2x_2}{x_0} - \frac{2y_2}{y_0} 
+ \frac{\beta-1}{2} \left( \beta- \frac{2x_1}{x_0} \right) \right)
.
\end{equation}

As before, the proof is divided in steps. 
In step 1 we make a case distinction regarding the sign of $z_1$, and characterise occurrence of the lcc in each case.
The lower and upper bound of the convergence exponent in the lcc is settles in step 2 and step 3, respectively.
In the last step, we finish the proof by showing how this relates to the actual statement of this theorem.

\subsubsection*{\underline{Step 1:} Growth of solutions}
We start with the difference equation, 
\begin{equation} \label{PP3}
\rho_{n+1} u_{n+2} + q_{n+1} u_{n+1} + \rho_n u_n = 0
.
\end{equation}
Proceeding as in the proof of \thref{P1} is not possible here, since we are in the case that the characteristic polynomial has a double zero.
Instead, set $r_i:=\frac{- q_i}{2 \rho_i }$
and divide \eqref{PP3} by $\rho_{n+1} \prod_{i=1}^{n+1} r_i$ to get
\begin{equation*} 
	\frac{u_{n+2} }{\prod_{i=1}^{n+1} r_i} - 2 \frac{u_{n+1}}{\prod_{i=1}^{n} r_i} + 
	\frac{\rho_n}{\rho_{n+1} r_n r_{n+1}}\frac{u_{n} }{\prod_{i=1}^{n-1} r_i} =0
	.
\end{equation*}
Introducing the new variable $v_n:= u_n \big/ \prod_{i=1}^{n-1} r_i$ and 
setting $C_n:=1-\frac{\rho_n}{\rho_{n+1} r_n r_{n+1}}$ gives
\begin{equation} \label{P4}
	v_{n+2} - 2 v_{n+1} + (1-C_n) v_{n} =0
	.
\end{equation}
A computation shows
\begin{equation} \label{P6}
	C_n=\frac{-z_1}{x_0 n}  + \frac{d}{n^2} + {\rm O}\Big(\frac{1}{n^{2+\epsilon}}\Big)
	,
\end{equation}
with some constant $d\in \bb R$. 

\paragraph*{Case 1: $z_1<0$.}
In this case, $J$ is in the lpc by Wouk's theorem, cf.\ \eqref{P2}.
\paragraph*{Case 2: $z_1>0.$}
Here, $\lim_{n\to \infty} n C_n = -z_1/x_0$ is negative
and \cite[Theorem~1,1.]{kooman:2007} gives two linearly independent solutions of \eqref{P4}, denoted by $(v_n^{(j)})_{n=1}^\infty$ for $j=1,2$, such that 
\begin{equation*} 
v_n^{(1)}=\overline{v_n^{(2)}} = (1+\mathrm{o}(1)) n^{1/4} \prod_{k=1}^{n-1} \big(1 + i \sqrt{-C_k} \big).  
\end{equation*}
The square of the absolute value of each factor is equal to
\begin{equation*}
	\left| 1 + i \sqrt{-C_k} \right|^2 =  1  - C_k 
	=  1 + 	\frac{z_1}{x_0 k}   + {\rm O}\Big(\frac{1}{k^{2}}\Big)  
	,
\end{equation*}
which leads to
\begin{equation*}
	\left| \prod_{k=1}^{n-1} \big(1 + i \sqrt{-C_k}\big) \right| 
	=  (c_1+o(1)) n^{\frac{z_1}{2x_0}}  
	= (c_1+o(1)) n^{ \frac{ x_1}{x_0} - \frac{ y_1}{y_0}- \frac{\beta}{2}} 
	,
\end{equation*} 
for some $c_1>0$, due to \cite[Lemma~4]{kooman:2007} adding a summable perturbation.
Thus, we get
\begin{equation} 
	\big| v_n^{(1)} \big|=\big| v_n^{(2)} \big| 
	=(c_1+o(1)) n^{ \frac{1}{4} + \frac{ x_1}{x_0} - \frac{ y_1}{y_0}- \frac{\beta}{2}}   \label{PP9}
		.  
\end{equation}

Substituting back via $u_n = v_n \prod_{i=1}^{n-1} r_i$ produces two solutions of \eqref{PP3},
denoted by $(u_n^{(j)})_{n=1}^\infty$ for $j=1,2$.
Again by \cite[Lemma~4]{kooman:2007}, we have
\begin{equation}\label{PP10}
	\bigg| 	\prod_{k=1}^{n-1} r_k   \bigg| = 
	\prod_{k=1}^{n-1} \Big(1 + \frac 1 k \Big(\frac{y_1}{y_0} - \frac{x_1}{x_0}\Big) + {\rm O}\Big(\frac{1}{k^{2}}\Big)   \Big)
	=(c_2+o(1)) n^{\frac{y_1}{y_0} - \frac{x_1}{x_0}}
	,
\end{equation}
for some $c_2> 0$. Together with \eqref{PP9} this results in the asymptotic behaviour
\begin{equation*} 
	\big| u_n^{(1)}  \big| = 	\big| u_n^{(2)}  \big| =  
	\big| v_n^{(1)}   \big| \bigg| \prod_{k=1}^{n-1} r_k  \bigg| =  (c_3+o(1))n^{\frac{1}{4}-\frac{\beta}{2}}
	,
\end{equation*}
where $c_3=c_1 c_2 >0$. In particular, $J$ is in the lcc if and only if $\beta>\frac{3}{2}$.

\paragraph*{Case 3: $z_1=0.$}
In that case, we have $\lim_{n\to \infty} n^2 C_n = d$, cf.\ \eqref{P6}. A calculation shows
\[
	d=\frac{2y_2}{y_0} - \frac{2x_2}{x_0} + (\beta-1)\frac{y_1}{y_0} + \frac{\beta ( \beta - 2)}{4}
	 = \frac{-z_2}{x_0} +  \frac{\beta ( \beta - 2)}{4}
	,
\]
where $z_2$ is defined in \eqref{P26}.
Note that $\beta\leq 2$ already implies that $J$ is in the lpc by Wouk's theorem since $z_1=0$, cf.\ \eqref{P2}.
We denote by $\alpha_1$ and $\alpha_2$ the zeros of the equation $X^2-X-d=0$, i.e.,
\[
	\alpha_{1,2}
	= \big(1 \pm \sqrt{1+4d}\,\big)\big/2
	.
\]
For $\alpha_1 \neq \alpha_2$, there are two linearly independent solutions of \eqref{P4} such that
\[
	v_n^{(1)} = (1+{\rm o}(1)) n^{\alpha_1}
	, \quad
	v_n^{(2)} = (1+{\rm o}(1)) n^{\alpha_2}
	.
\]
This follows from either \cite[Theorem~10.1,(1)]{kooman:1998}, or \cite[Theorem~1,2.]{kooman:2007}.
Actually, the case $d=0$ is already treated in \cite[Theorem~10.3]{coffman:1964}.

In the case of a double zero $\alpha_1=\alpha_2=1/2$, we get two solutions of \eqref{P4} with
\[
	v_n^{(1)} = (1+{\rm o}(1)) n^{1/2}, \quad
	v_n^{(2)} = (1+{\rm o}(1)) \log(n) n^{1/2}.
\]
To transform these solutions back 
to solutions $u_n^{(j)}$ of \eqref{PP3},
note that 
\[
\bigg| 	\prod_{k=1}^{n-1} r_k  \bigg| = (c_3+{\rm o}(1))  n^{-\beta/2}
,
\]
by \eqref{PP10} together with $z_1=0$. 
\begin{enumerate}[\text{Case} 3a:]
	\item $d<-1/4.$ \ \
	In this case, $\alpha_1$ and $\alpha_2$ 
	are two distinct complex conjugate numbers with $\Re \alpha_i = 1/2$, and we get
	two solutions of \eqref{PP3} with
	\[
		|u_n^{(1)}| = |u_n^{(2)}|=(1+{\rm o}(1)) n^{(1-\beta)/2}
		.
	\]
	Thus, $J$ is in the lcc if and only if $\beta>2$.
	\item $d=-1/4. $ \ \
	Here, $\alpha_1=\alpha_2=1/2$ is a double zero, and we get
	\[
		|u_n^{(1)}| = (1+{\rm o}(1))  n^{(1-\beta)/2}\log n, \quad
		|u_n^{(2)}| = (1+{\rm o}(1)) n^{(1-\beta)/2}
		.
	\]	
	As before, $J$ is in the lcc if and only if $\beta>2$.
	\item $d>-1/4.$ \ \ 
	In that case, $\alpha_1$ and $\alpha_2$ are two distinct real zeros, and we get two solutions of \eqref{PP3} such that
	\[
		|u_n^{(1)}| = (1+{\rm o}(1)) n^{(1 + \sqrt{1+4d} - \beta)/2}, \quad
		|u_n^{(2)}| = (1+{\rm o}(1)) n^{(1 - \sqrt{1+4d} - \beta)/2}.
	\]
	Here, $J$ is in the lcc if and only if the dominating solution is square-summable, i.e,
	\[
		 \sqrt{1+4d}  < \beta -2
		.
	\]
	For $\beta\leq 2$ this inequality is obviously false, i.e., $J$ is in the lpc.
	If $\beta>2$, then the above condition is further equivalent to 
	$\beta < \frac{3}{2} + \frac{2z_2}{x_0}$.
\end{enumerate}

%
%

\subsubsection*{\underline{Step 2:} The lower bound in the lcc}

This step can be done exactly as in \thref{P1}.
When $J$ is in the lcc, we get as before  $\rho(\sigma(J))\geq 1 /\beta$, as well as
\[
	\limsup_{r\to \infty} \frac{n_{\sigma} (r)}{r^{1/{\beta}}}
	\geq  \frac{\beta - 1}{\beta} 
	\Big(\frac{1}{x_0}\Big)^{1/\beta}
	.
\]

\subsubsection*{\underline{Step 3:} The upper bound in the lcc}

In the first step we have seen that the difference equation \eqref{PP3} has a fundamental solution
$u_n^{(1)},u_n^{(2)}$ such that the dominating solution satisfies
$|u_n^{(1)}| \asymp \lambda(n)$ where either $\lambda(n):=n^\gamma$ or $\lambda(n):=n^\gamma \log n$ for some $\gamma \in \bb R$.


Recall that the orthogonal polynomials of the first and second kind, denoted by $P_n(0)$ and $Q_n(0)$, respectively, are linearly independent solutions of \eqref{P3}. 
The quotient $(\left|P_n(0)\right| + \left|Q_n(0)\right|)/{\lambda(n)}$
is bounded from above since $P_n(0)$ and $Q_n(0)$
can be written as linear combinations of $u_n^{(1)}$ and $u_n^{(2)}$. 
It is also bounded away from zero, since $u_n^{(1)}$ is
a linear combination of $P_n(0)$ and $Q_n(0)$ and 
$ | u_n^{(1)}|	/{\lambda(n)}$ is bounded away from zero.
Thus, we obtain $P_n(0)^2  + Q_n(0)^2  \asymp \lambda(n)^2$.

Now we write the Jacobi matrix as a Hamburger Hamiltonian of a canonical system, cf.\ \cite{pruckner.romanov.woracek:2016} or \cite{kac:1999} for details about this reformulation.
We denote by $(l_n)_{n=1}^\infty$ and $(\phi_n)_{n=1}^\infty$ the sequences of lengths and angles of the corresponding Hamburger Hamiltonian.
By \cite[(1.5),(1.6)]{pruckner.romanov.woracek:2016}, we have that
\begin{align*}
	&\, l_n= P_n(0)^2  + Q_n(0)^2  \asymp \lambda(n)^2 
	, \\
	&\big|\sin(\phi_{n+1}- \phi_n ) \big| = 
	1\big/\big(  \rho_n  \sqrt{l_n l_{n+1}}\big)
	\asymp n^{-\beta} \lambda(n)^{-2}	
	.
\end{align*}
With the notation from \cite{pruckner.romanov.woracek:2016}, 
the lengths and angle-differences are regularly distributed. 
Moreover, we have $\Delta_l=-2\gamma$ and $\Delta_\phi=\beta+2\gamma$, both expressions exist as a limit and 
$\Delta_l + \Delta_\phi = \beta$.

\paragraph*{Case 2: $z_1>0.$}
In this case we have $\gamma = (1 -2\beta)/4$.  
For $3/2 < \beta < 2$, we have $\Delta_l + \Delta_\phi = \beta <2$.
By \cite[Theorem~2.7]{pruckner.romanov.woracek:2016} 
the order of $B$, i.e.\ $\rho(\sigma(J))$, does not exceed
\[
	\frac{1-\Delta_\phi - \frac{\Lambda}{2}}{\Delta_l - \Delta_\phi + \Lambda} =\frac{1}{2(\beta - 1)}
	.
\]
For $\beta\geq 2$, \cite[Theorem~2.22,(i)]{pruckner.romanov.woracek:2016}
is applicable and gives $\rho(\sigma(J))=1/\beta_1$.

\paragraph*{Case 3: $z_1=0.$}
In this case, $\beta>2$ is necessary for occurrence of the lcc by Wouk's theorem.
Due to $\Delta_l + \Delta_\phi = \beta >2$, \cite[Theorem~2.22,(i)]{pruckner.romanov.woracek:2016} is applicable and gives $\rho(\sigma(J))=1/\beta_1$.



\subsubsection*{\underline{Step 4:} Conclusion}
Fist consider $\beta >\frac{2x_1}{x_0} -\frac{2y_1}{y_0}$ which is equivalent to $z_1<0$.
Hence, we are in case 1, and $J$ is in the lpc by the first step.

Similarly,  $\beta <\frac{2x_1}{x_0} -\frac{2y_1}{y_0}$ is equivalent to $z_1>0$, which is case 2. 
By the first step, $J$ is in the lcc if and only if $\beta > 3/2$. Regarding the convergence exponent, \eqref{P27} holds by the third step.

\noindent
Therefore, case $(ii)$ in the formulation of the theorem is settled, as well as case $(i)$ with the possible exception of $\beta=\frac{2x_1}{x_0} -\frac{2y_1}{y_0}\leq \frac 3 2$.
In this case we have $z_1=0$, i.e.\ we are in case 3. Due to $\beta \leq \frac 3 2 $, $J$ is in the lpc by Wouk's theorem.

Thus, it remains to treat case $(iii)$, i.e. $\beta=\frac{2x_1}{x_0} -\frac{2y_1}{y_0}> \frac 3 2$.
Once more we have $z_1=0$ and, thus, fall into case 3.
Recall from the first step that $J$ is in the lcc if and only if 
\begin{equation*} 
	d\leq -\tfrac{1}{4} \,\text{ and }\, \beta > 2,
	\quad \text{ or, } \quad
	 2<\beta< \tfrac{3}{2} + \tfrac{2 z_2}{x_0}
	.
\end{equation*}
Next, we show that $d \leq -1/4$ and $\beta>2$ already implies $\beta< \frac{3}{2} + \frac{2 z_2}{x_0}$. 
By solving a quadratic equation one can show that $d = -\tfrac{z_2}{x_0} + \tfrac{\beta (\beta-2)}{4} \leq -\frac{1}{4}$
implies
\[
	\beta \leq 1 + 2 \big(\tfrac{z_2}{x_0} \big)^{1/2}
	.
\]
For $z_2/x_0 = 1/4$ this would give $\beta \leq 2$, which would contradict $\beta >2$.
Thus, we have $z_2/x_0 \neq 1/4$ and obatin, again by solving a quadratic equation, the estimate
\[
	\beta\leq 1 + 2 \big(\tfrac{z_2}{x_0} \big)^{1/2} <  \tfrac{3}{2} + \tfrac{2 z_2}{x_0}
	.
\]
Hence, occurrence of the lcc in case 3 is equivalent to $2< \beta < \tfrac{3}{2} + \tfrac{2 z_2}{x_0}$.
In the lcc, we have $\rho(\sigma(J))=1/\beta$ by the third step.
\qed

\smallskip

\begin{remark} \thlab{P28}
The techniques used in the third step of the proof of \thref{P1} does not seem to be suitable in the situation of \thref{ADD1}.

To demonstrate this, consider the case 2, i.e.\ $z_1>0$. By the first step, we know $P_n^2,Q_n^2 = {\rm O }(n^{(1-2\beta)/2})$, and \cite[Theorem~1.2]{berg.szwarc:2014} gives $\rho(\sigma(J))\leq 1/(\beta-\frac{1}{2})$.
Together with the lower estimate, we get $\rho(\sigma(J))\in [1/\beta, 1/(\beta-\frac{1}{2})]$.
Hence, contrary to the situation in \thref{P1}, this only shows that the convergence exponent is contained in an interval.
Also the estimate of the Nevanlinna matrix, as performed in the proof of \thref{P1}, does not improve the size of this interval.

Using \cite[Theorem~2.7]{pruckner.romanov.woracek:2016} improves our result drastically:
For $\beta<2$ the size of the interval shrinks, and for $\beta \geq 2$ this determines the convergence exponent.
\end{remark}

\section*{Acknowledgement}
I thank the referee for the constructive comments and suggestions which improved the result and shortened the proof.

Special thanks go to R.Romanov for helpful discussions about Theorem 1.

\medskip
\noindent
This work was supported by the Austrian Science Fund [FWF, I\,1536--N25; P\,30715-N35].





\ifthenelse{\Draft=1}{
\newpage

\noindent
{\large\sc labels:}\\[10mm]
	
P:\hspace*{3mm}
\framebox{
\begin{tabular}{r@{\ }r@{\ }r@{\ }r@{\ }r@{\ }r@{\ }r@{\ }r@{\ }r@{\ }r@{\ }}
	** &  1 &  2 &  3 &  4 &  5 &  6 &  7 & 8 &  9 \\
	 10 &  11 & 12 &  13 &  14 & 15 &  16 & 17  &  18 & 19 \\
	 20 &  21 & 22 &  23 &  24 &  25 &  26 &  27 &  28 &  . \\
	 . &  32 &  . &  . &  . &  . &  . &  . &  . &  . \\
\end{tabular}
}
\\[1cm]
}
{}

\end{document}